\documentclass[12 pt]{article}

\usepackage{fancyhdr}
\usepackage[noend]{algorithmic}
\usepackage{algorithm}
\usepackage{amscd}
\usepackage{amsfonts}
\usepackage{amsmath}
\usepackage{amssymb}
\usepackage{amsthm}
	\newtheorem{theorem}{Theorem}
	\newtheorem{definition}{Definiton}
\usepackage{blindtext}
\usepackage[font = {small}, justification = centering]{caption}
\usepackage{grffile}
\usepackage{mathtools}
\usepackage{multicol}
\usepackage{pgf}
\usepackage{tikz}
\usepackage[utf8]{inputenc}
\usepackage{geometry}
\usepackage{caption}
\usepackage{subcaption}

\usetikzlibrary{arrows}
\usetikzlibrary{backgrounds}
\usetikzlibrary{calc}
\usetikzlibrary{decorations.pathmorphing}
\usetikzlibrary{fit}
\usetikzlibrary{petri}
\usetikzlibrary{positioning}
\usetikzlibrary{trees}

\graphicspath{ {} }
\DeclareGraphicsExtensions{.jpeg,.jpg,.png,.pdf}

\DeclareMathOperator{\N}{\mathbb{N}}

\DeclareMathOperator{\Z}{\mathbb{Z}}

\title{A Novel AQC Factoring Algorithm}
\author{Matthew B. Crawford\thanks{Naval Surface Warfare Center Dahlgren Division  (NSWCDD), matthew.b.crawford1@navy.mil}}
\date{}

\newcount\mycount

\begin{document}
\maketitle

\makeatletter

%
%
%


\begin{abstract}
Due to recent technological advances, actual quantum devices are being constructed and used to perform computations.  As a result, many classical problems  are being restated so as to be solved on quantum computers.  Some examples include satisfiability problems \cite{farhi2001quantum}; clustering and classification, \cite{lloyd2013quantum}, \cite{lloyd2014quantum}, \cite{rebentrost2014quantum}; protein folding \cite{perdomo2012finding}; and simulating many-body systems \cite{bloch2012quantum}.  Converting these classical problems to a quantum framework is not always straightforward.  As such, instances where researchers explicitly elucidate the conversion process are not only valuable in their own right, but are likely to spawn  new ideas and creative ways in regards to problem solving.  In this paper, we propose a classical factoring algorithm, which we then convert into a quantum framework.  Along the way, we discuss the subtle similarities and differences between the approaches, and provide a general comparison of their performance. It is our desire to not only introduce an interesting approach to factoring, but to hopefully promote more creative ways to solving problems using quantum computers.  The key to our algorithm is that we convert the factoring problem to a graph theory problem using elements from group theory.  The move to a graph-theoretic approach ultimately eases the transition to a quantum setting.\\

\end{abstract}

\newpage

\tableofcontents

\newpage

\listoffigures


\listoftables


\listofalgorithms

\newpage

\section{Introduction}

Due to recent technological advances, actual quantum devices are being constructed and used to perform computations.  As a result, many classical problems  are being restated so as to be solved on quantum computers.  Some examples include satisfiability problems \cite{farhi2001quantum}; clustering and classification, \cite{lloyd2013quantum}, \cite{lloyd2014quantum}, \cite{rebentrost2014quantum}; protein folding \cite{perdomo2012finding}; and simulating many-body systems \cite{bloch2012quantum}.  Converting these classical problems to a quantum framework is not always straightforward.  As such, instances where researchers explicitly elucidate the conversion process are not only valuable in their own right, but are likely to spawn  new ideas and creative ways in regards to problem solving.  In this paper, we propose a classical factoring algorithm, which we then convert into a quantum framework.  Along the way, we discuss the subtle similarities and differences between the approaches, and provide a general comparison of their performance. It is our desire to not only introduce an interesting approach to factoring, but to hopefully promote more creative ways to solving problems using quantum computers.  The key to our algorithm is that we convert the factoring problem to a graph theory problem using elements from group theory.  The move to a graph-theoretic approach ultimately eases the transition to a quantum setting.\\
 The remainder of the paper is organized as follows.  In section 2, we present a graph-theoretic approach to factoring integers.  Given a semiprime $N,$ we build a graph with $N-1$ vertices representing the integers $1,\ldots ,N-1.$  We develop an equivalence relation on $1,\ldots,N-1$ that partitions the graph into disjoint connected components.  By carefully constructing a graph in this way, we show that one can (with high probability) find a component that contains one of the prime factors of $N.$  The constructions in section 2 will be used in sections 3 and 4, where we discuss a random walk and an adiabatic quantum (AQC) approach to selecting the correct components.  Section 5 presents a series of experiments where we test our algorithms and compare their performance.  We conclude in section 6, where we discuss plans for future work.


\section{Number and Graph Theory}
We begin by defining a family of graphs derived from a semiprime.

\begin{definition}\label{def:graph}
Let $N=pq$ for two distinct primes $p$ and $q$, and let $\alpha \in (\Z/N\Z)^*$. We define $G_{N,\alpha}$ to be the graph on vertices labeled $1,\ldots,N-1$ with edges between $i$ and $j$ if and only if $j=\alpha i \mod N$.
\end{definition}  For example, given $N=15$, and $\alpha=2$, we have Figure \ref{fig:2.1} below.

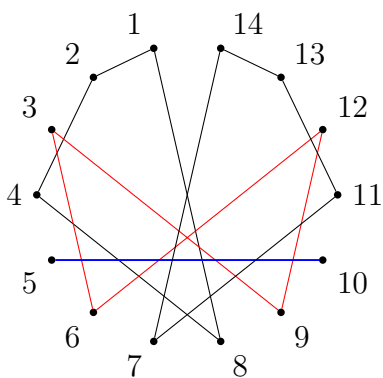
\begin{figure}[h!]
\centering
\begin{tikzpicture}
\def \margin {4} 

\foreach \s in {1,...,14} 
	{
	\node[circle,fill=black,inner sep = 1 pt, minimum size=0 pt, outer sep=0 pt, label=360/14 *((\s-1/2)+14/4):$\s$] (\s) at ({360/14 *((\s-1/2)+14/4)}:2) {};
  	}
\draw[black] (1) -- (2) -- (4) -- (8) -- (1);
\draw[red] (3) -- (6) -- (12) -- (9) -- (3);
\draw[blue] (5) -- (10) -- (5);
\draw[black] (7) -- (14) -- (13) -- (11) -- (7);
\end{tikzpicture}
\caption{Graph for $N=15$ and $\alpha=2$.}
\label{fig:2.1}
\end{figure}

\begin{figure}[h!]
\centering
\begin{tikzpicture}
\def \margin {4} 

\foreach \s in {1,...,54}
	{
\node[circle,fill=black,inner sep = 1 pt, minimum size=0 pt, outer sep=0 pt, label=360/54*((\s-1/2)+54/4):$\s$] (\s) at ({360/54*((\s-1/2)+54/4)}:4) {};
	}
\draw[black] (1) -- (3) -- (9) -- (27) -- (26) -- (23) -- (14) -- (42) -- (16) -- (48) -- (34) -- (47) -- (31) -- (38) -- (4) -- (12) -- (36) -- (53) -- (49) -- (37) -- (1);
\draw[black] (2) -- (6) -- (18) -- (54) -- (52) -- (46) -- (28) -- (29) -- (32) -- (41) -- (13) -- (39) -- (7) -- (21) -- (8) -- (24) -- (17) -- (51) -- (43) -- (19) -- (2);
\draw[red] (5) -- (15) -- (45) -- (25) -- (20) -- (5);
\draw[red] (10) -- (30) -- (35) -- (50) -- (40) -- (10);
\draw[blue] (11) -- (33) -- (44) -- (22) -- (11);
\end{tikzpicture}
\caption{Graph for $N=55$ and $\alpha=3$.}
\label{fig:2.2}
\end{figure}
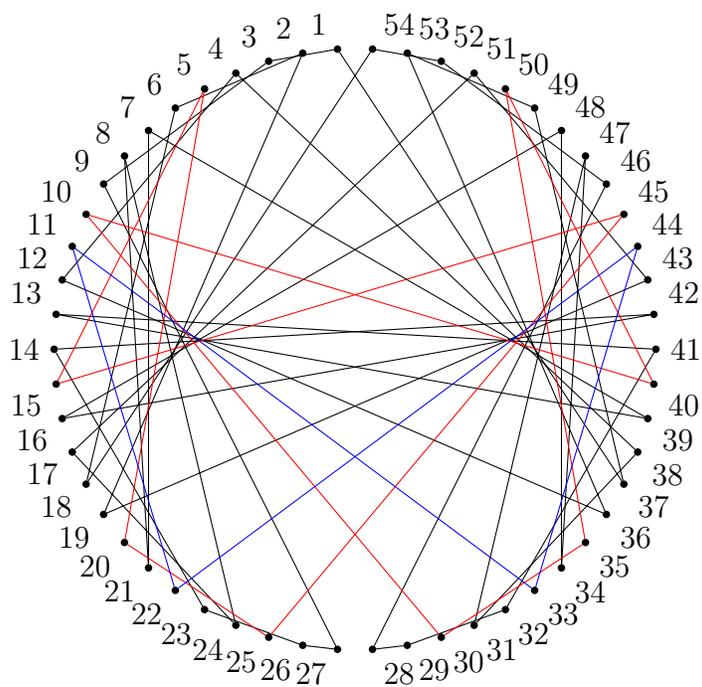

We can see that there are cycles of different types; the cycles containing multiples of one of the factors of $N$ (red), the cycles containing multiples of the other factor of $N$ (blue), and the cycles containing elements relatively prime to $N$ (black).\\ \\

\begin{theorem}\label{theorem:cycles}
Let $G=G_{N,\alpha}$ be the graph defined in Definition \ref{def:graph} for some distinct primes $p$ and $q$ and some totative $\alpha $ of $N$. Let $\Z/N\Z(x)_\alpha$ denote the orbit of $x \in \Z/N\Z$ under the action of multiplication by $\alpha,$ and let $[a,b]$ denote the Least Common Multiple of $a$ and $b$. Then the number of cycles $C_G$ in $G$  is:
\begin{equation}
C_G=\frac{q-1}{|\Z/N\Z(p)_\alpha|}+\frac{p-1}{|\Z/N\Z(q)_\alpha|}+\frac{(p-1)(q-1)}{[|\Z/N\Z(p)_\alpha|,|\Z/N\Z(q)_\alpha|]}.
\end{equation}
\end{theorem}

\begin{proof}
The first summand is the number of cycles containing vertices $i$ with $\gcd(i,N)=p$. These $i$ form a subgraph of $G$ where $\gcd(\frac{i}{p},\frac{N}{p})=\gcd(\frac{i}{p},q)=1$. Thus these vertices corresponding to elements of $\Z/N\Z$ biject with the elements of the group $(\Z/q\Z)^*$ by dividing by $p$. By the Orbit-Stabilizer Theorem, letting our group be $(\Z/q\Z)^*$, and $(\Z/q\Z)^*_x$ be the stabilizer of $x$ in $(\Z/q\Z)^*$, we have:
\begin{align*}
|(\Z/q\Z)^*|&=|(\Z/q\Z)^*_1||(\Z/q\Z)^*(1)_\alpha|\\
\implies |(\Z/q\Z)^*_1|&=\frac{q-1}{|(\Z/q\Z)^*(1)_\alpha|}\\
&=\frac{q-1}{|\Z/N\Z(p)_\alpha|}.
\end{align*}
Because the order of the group is the number of vertices, $i$,  such that $\gcd(i,N)=p$, and the size of the orbit is the number of aforementioned vertices in a cycle with $p$, and since each orbit has the same size, the stabilizer must be the number of cycles.\\ \\
We recognize the second summand in the same way, but with vertices $i$ such that $\gcd(i,N)=q$, and similarly $|(\Z/p\Z)^*_1|=\frac{p-1}{|\Z/N\Z(q)_\alpha|}$.\\ \\
The last summand is the set of vertices $i$ such that $\gcd(i,N)=1$. In the same way as the other two summands, we use the Orbit-Stabilizer Theorem with $(\Z/N\Z)^*$ and arrive at:
\begin{align*}
|(\Z/N\Z)^*_1|&=\frac{(p-1)(q-1)}{|(\Z/N\Z)^*(1)_\alpha|}\\
&=\frac{(p-1)(q-1)}{|(\Z/p\Z)^*\times(\Z/q\Z)^*(1,1)_\alpha|}\\
&=\frac{(p-1)(q-1)}{\left [ |(\Z/p\Z)^*(1)_\alpha|,|(\Z/q\Z)^*(1)_\alpha| \right ]}\\
&=\frac{(p-1)(q-1)}{[|\Z/N\Z(p)_\alpha|,|\Z/N\Z(q)_\alpha|]}.
\end{align*}
Thus the total number of cycles is as stated above.
\end{proof}



Given a graph constructed in this manner, we may randomly select, or \emph{mark}, $k \in \N $ vertices, and check whether the following conditions are satisfied:
\begin{table}[h!]
\hspace*{-2 mm}
\begin{tabular}{|l|}
	\hline
	\textbf{Condition 1} At least one vertex in each black cycle is marked,\\
	\textbf{Condition 2} At least one blue or red cycle contains no marked vertices.\\
	\hline
\end{tabular}
\caption{Conditions to be satisfied when marking graph vertices.}
\label{table:conditions}
\end{table}

 Let $P(k)$ denote the probability of marking $k$ vertices chosen according to the uniform distribution on $G_{N,\alpha}$ such that the conditions above are satisfied. As an example, consider Figure \ref{fig:2.1} (See Table \ref{table:2.3} for the corresponding probabilities).

\begin{table}[h!]
\hspace*{-2 mm}
\begin{tabular}{|c||c|c|c|c|c|c|c|c|c|c|c|c|c|c|}
\hline
$k$ & 2 & 3 & 4 & 5 & 6 & 7 & 8 & 9 & 10 & 11 & 12\\
\hline
\hline
$P(k)$ & $.176$ & $.396$ & $.468$ & $.432$ & $.349$ & $.259$ & $.179$ & $.115$ & $.067$ & $.033$ & $.011$ \\
\hline
\end{tabular}
\caption{Probabilitites $P(k)$ of the $k$ marked vertices satisfying the conditions.}
\label{table:2.3}
\end{table}

\begin{figure}[h!]
\begin{center}
\includegraphics[width=9 cm, height= 6.75 cm]{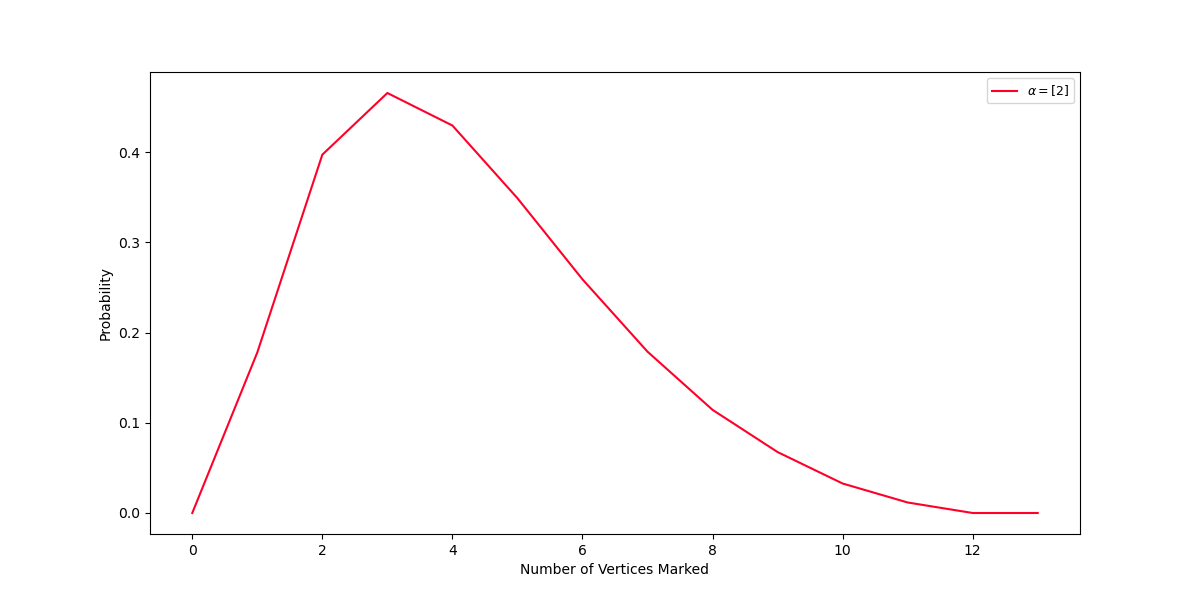}
\end{center}
\caption{Probabilities of the marked vertices meeting the conditions in Table \ref{table:conditions} for $G_{15,2}$.}
\label{fig:2.4}
\end{figure}

In this case, we exclude $k=1$,13, and 14.  For $k=1$, there are two black cycles, so it is impossible to mark both with one vertex.  For $k=13$ and 14, there are twelve vertices total in both black and red cycles, so it is impossible to mark another non-blue vertex.\\ \\
We may also use multiple $\alpha$'s, thus possibly connecting disjoint cycles. However, it is important to notice that only cycles of the same color can be connected by using multiple $\alpha$'s.  Figure \ref{fig:2.5} shows the graph of $N=15$ with $\alpha=2$ and $\alpha=7,$ and Table \ref{table:2.6} contains the corresponding probabilities.  We note that now $P(1)$ is nonzero because there is only one black cycle. We also note that the maximum probability using both values of $\alpha$ is greater than the probability of only using $\alpha=2$. This is because using $\alpha=2$ and $\alpha=7$ connects the original two black cycles into one component, and leaves the other two red and black cycles unchanged. Thus, the probability of any $k$ random vertices meeting the conditions increases.  We will show later that the conditions in Table \ref{table:conditions} are too strong; i.e., we may ``succeed'', even without the conditions being met.


\begin{figure}[h!]
\centering
\begin{tikzpicture}
\def \margin {4} 

\foreach \s in {1,...,14} 
	{
	\node[circle,fill=black,inner sep = 1 pt, minimum size=0 pt, outer sep=0 pt, label=360/14 *((\s-1/2)+14/4):$\s$] (\s) at ({360/14 *((\s-1/2)+14/4)}:2) {};
  	}
\draw[black] (1) -- (2) -- (4) -- (8) -- (1);
\draw[red] (3) -- (6) -- (12) -- (9) -- (3);
\draw[blue] (5) -- (10) -- (5);
\draw[black] (7) -- (14) -- (13) -- (11) -- (7);
\draw[red] (3) -- (6) -- (12) -- (9) -- (3);
\draw[black] (2) -- (14) -- (8) -- (11) -- (2);
\draw[blue] (10) -- (10);
\draw[black] (1) -- (7) -- (4) -- (13) -- (1);
\draw[blue] (5) -- (5);
\end{tikzpicture}
\caption{Graph for $N=15$ and $\alpha=2,7$.}
\label{fig:2.5}
\end{figure}

\begin{figure}[h!]
\begin{center}
\includegraphics[width=9 cm, height=6.75 cm]{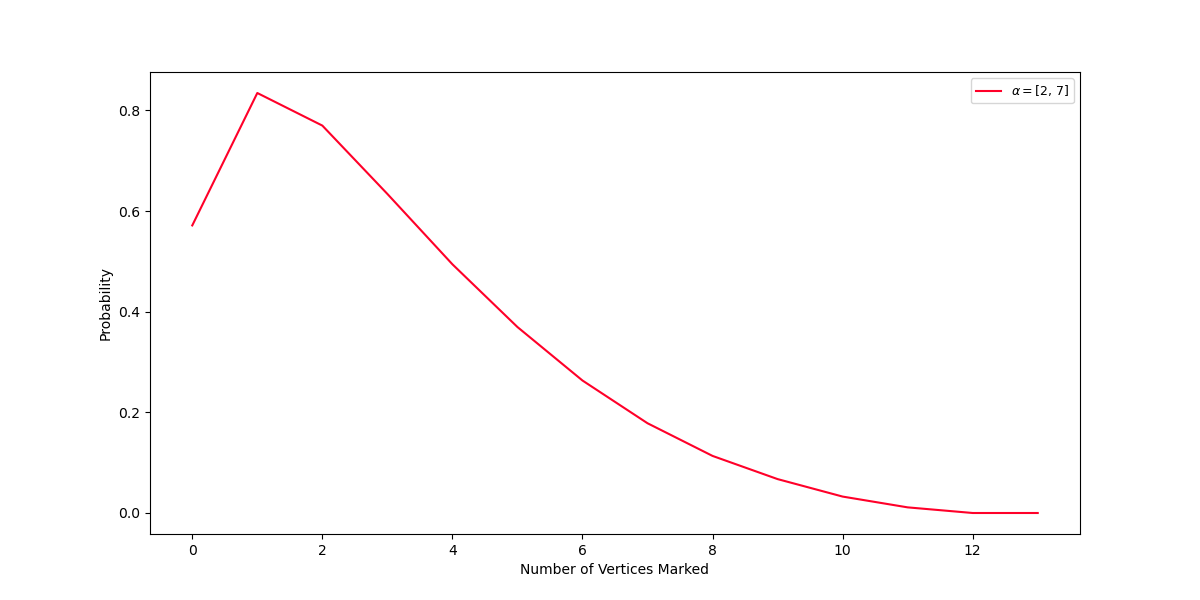}
\end{center}
\caption{Probabilities of the marked vertices meeting the conditions in Table \ref{table:conditions} for $G_{15,2} \cup G_{15,7}$.}
\label{fig:2.7}
\end{figure}

\begin{table}[h!]
\hspace*{-9 mm}
\begin{tabular}{|c||c|c|c|c|c|c|c|c|c|c|c|c|c|c|}
\hline
$k$ & 1 & 2 & 3 & 4 & 5 & 6 & 7 & 8 & 9 & 10 & 11 & 12\\
\hline
\hline
$P(k)$ & $.581$ & $.835$ & $.768$ & $.632$ & $.494$ & $.366$ & $.266$ & $.180$ & $.115$ & $.067$ & $.033$ & $.011$ \\
\hline
\end{tabular}
\caption{Probabilitites $P(k)$ of the $k$ marked vertices satisfying the conditions.}
\label{table:2.6}
\end{table}


\section{Classical Approach}

In this section, we describe a random walk approach for determining the probabilities discussed in the previous section.  We begin by considering an arbitrarily large number of walkers placed uniformly on the graph $G_{N,\alpha}$ for some $N$ and $\alpha$. We then allow the walkers to move about randomly along connected components.
  Let $M$ be a random set of $1 \leq k \leq N-1$ vertices to mark.  For each vertex $m \in M$, replace $m$ with a {\it wormhole} with the following property:  any walker that enters the wormhole is  instantaneously teleported to a new vertex which is not a wormhole. We define the probability of a walker being teleported to any particular vertex as the proportion of walkers currently occupying that vertex (See Theorem \ref{theorem:graphmen probs} below). This operation of changing vertices into wormholes changes the distribution of walkers on $G_{N,\alpha}$.  As a result, walkers will tend to accumulate on unmarked components.  By carefully choosing the vertices to mark, we can ensure (with high probability) that walkers accumulate on components that contain a multiple of one of the factors of $N.$
\par As an example, consider $G_{15,2}$ in Figure \ref{fig:2.1}.  If $M=\{1,14\}$, over time any walker who walks on to vertex or $1$ or $14$ would be teleported to a vertex in a red or blue cycle.  This walker would then be stuck there, because none of the vertices in the red and blue cycles are wormholes.  Over time, the distribution of walkers will tend to concentrate on unmarked components.

\par Here we provide an example that demonstrates why our earlier conditions for ``success'' were too strong (See Table \ref{table:conditions}).  Recall that, given some collection of vertices $M,$ our goal is to ultimately have the distribution of walkers solely on red or blue cycles.  We will show that there exist collections of vertices that do not meet the conditions, but that still result in the distribution of walkers to be on only red or blue cycles.  Considering $G_{15,2}$ again, if $M=\{1,2,3,4,5,7,8,11,13,14\}$, then even though there are wormholes in every cycle, the only vertices that are not wormholes are $6,9,10$, and $12$ -- vertices in red or blue cycles. Therefore, over time, all of the walkers will be transported to these vertices.\\ \\

\begin{theorem}\label{theorem:graphmen probs}
Let $N=pq$, where $p$ and $q$ are distinct primes. Let $\alpha$ be a totative of $N$. Let $G_{N,\alpha}$ be the graph defined in Definition \ref{def:graph}. Let $M$ be a set of $1 \leq k \leq N-1$ vertices to mark on $G_{N,\alpha}$ as wormholes. Let $P(t)_{r \not\in M}$ be the probability distribution of walkers on the unmarked vertices of $G_{N,\alpha}$ at time $t$. Let the distribution of walkers be uniformly distributed at time $t=0$, that is, each entry of $P(0)_{r \not\in M}$ is equal to $\frac{1}{N-1-k}$. Let the walking matrix $E=E_{ij}$ be the probability of a walker on vertex $i$ going to vertex $j$ via an edge and not a wormhole. Then:
\begin{equation*}
P(t)_{r \not\in M}=\left(1+\frac{W}{S}\right)EP(t-1)_{r \not \in M},
\end{equation*}
where $S=||EP(t-1)_{r \not \in M}||_1$, and $W=1-S$.
\end{theorem}

\begin{proof}~\\
The travels of the walkers that do not go through wormholes is modeled by $EP(t-1)_{r \not \in M}$.  Therefore, $S$ is the proportion of walkers who do not go in a wormhole, while $W$ is the proportion of walkers who do. The $W$ walkers who do go into a wormhole are then redistributed to the other non-wormhole vertices based on their population proportion, which can be found by normalizing the distribution of walkers from edge-traveling only.
\end{proof}

\begin{theorem}\label{theorem:graphmen and hamiltonian}
Let $N=pq$, where $p$ and $q$ are distinct primes. Let $\alpha$ be a totative of $N$. Let $G_{N,\alpha}$ be the graph defined in Definition \ref{def:graph}. Let $M$ be a set of $1 \leq k \leq N-1$ vertices to mark on $G_{N,\alpha}$ as wormholes. 
Let $H_p$ be the combinatorial Laplacian for the graph $G_{N,\alpha},$
and let $P(t)_{r \not\in M}$ be the probability distribution of walkers on the unmarked vertices of $G_{N,\alpha}$ at time $t$. Then, if $P(0)_{r \not\ in M}$ is the uniform distribution and $H_p$ has a non-zero spectral gap, $$\lim_{t \to \infty} P(t)_{r \not\in M}=\frac{v_{\min}}{||v_{\min}||_1},$$ where $v_{\min}$ is eigenvector of the smallest eigenvalue of $H_p$. 
\end{theorem}

\begin{proof}~\\
Consider the differential equation below, where $H$ is a real hermitian matrix:
\begin{equation*}
\frac{df}{dt}=-Hf.
\end{equation*}
It is a well known fact that the general solution to this differential equation is:
\begin{equation*}
f(t)=\sum_{i=1}^{n}C_ie^{-\lambda_i t}v_i,
\end{equation*}
where the $\lambda_i$'s are eigenvalues of $H$, and the $v_i$'s are corresponding orthonormal eigenvectors. We may assume the terms are ordered in increasing order by eigenvalue. 
\par Now, a first order approximation of $f(t+\Delta t)$ is:
\begin{align*}
f(t+\Delta t) &\approx \frac{df}{dt}\Big|_{t}\Delta t + f(t)\\
&=-Hf(t)\Delta t + f(t)\\
&=(I-H\Delta t)f(t).
\end{align*}
As $t$ goes to infinity, $f(t)$ will converge to 0, as will this approximation. However, if we normalize  $f(t)$ (in $\ell_1$) after each time step, we will instead converge to $\frac{v_1}{||v_1||_1}$.  Let $t_0$ represent the initial time.  Then the above approximation becomes
$$f(t_0 + \Delta t) \approx (I - \Delta t H)f(t_0) = \sum_{i\ge 1} C_i(1 - \lambda_i \Delta t)e^{-\lambda_i t_0}v_i,$$
which, after normalization, becomes
$$\frac{f(t_0 + \Delta t)}{||f(t_0 + \Delta t)||_1}=\frac{\sum_{i \ge 1}C_i(1 - \lambda_i \Delta t)e^{-\lambda_i t_0}v_i}{\left | \left | \sum_{i \ge 1}C_i(1 - \lambda_i \Delta t)e^{-\lambda_i t_0}v_i \right | \right |_1}.$$
Proceeding to the next iteration gives
\begin{eqnarray*}
	f(t + 2 \Delta t) &=& (I - \Delta t)H \frac{f(t_0 + \Delta t)}{||f(t_0 + \Delta t)||_1}\\
&=&\frac{\sum_{i\ge 1}C_i (1 - \lambda_i \Delta )^2 e^{-\lambda_i t_0}v_i}{||\sum_{i\ge 1}C_i (1 - \lambda_i \Delta ) e^{-\lambda_i t_0}v_i||_1}\\
\end{eqnarray*}
Normalizing again gives
$$\frac{f(t + 2 \Delta t)}{\left | \left | f(t + 2 \Delta t) \right | \right |_1} = \frac{ \displaystyle \frac{\sum_{i\ge 1}C_i (1 - \lambda_i \Delta )^2 e^{-\lambda_i t_0}v_i}{||\sum_{i\ge 1}C_i (1 - \lambda_i \Delta ) e^{-\lambda_i t_0}v_i||_1} }{ \displaystyle \left | \left |\frac{\sum_{i\ge 1}C_i (1 - \lambda_i \Delta )^2 e^{-\lambda_i t_0}v_i}{||\sum_{i\ge 1}C_i (1 - \lambda_i \Delta ) e^{-\lambda_i t_0}v_i||_1 } \right | \right |_1 } = \frac{\sum_{i\ge 1}C_i (1 - \lambda_i \Delta )^2 e^{-\lambda_i t_0}v_i}{\left | \left | \sum_{i\ge 1}C_i (1 - \lambda_i \Delta )^2 e^{-\lambda_i t_0}v_i \right | \right |_1}.$$
Continuing inductively, we have the following for the $n^{th}$ iteration:
\begin{eqnarray*}
	\frac{f(t_0 + n\Delta t)}{||f(t_0 + n\Delta t)||_1}& = &\frac{\sum_{i\ge1} C_i(1-\lambda_i \Delta t)^n e^{-\lambda_i t_0}v_i}{\left | \left | \sum_{i\ge1} C_i(1-\lambda_i \Delta t)^n e^{-\lambda_i t_0}v_i \right | \right |_1 }\\
&=&\frac{v_1 + \sum_{i>1}\frac{C_i}{C_1} \left [ \frac{1 - \lambda_i \Delta t}{1 - \lambda_1 \Delta t}\right ]^n e^{(\lambda_1 - \lambda_i)t_0}v_i}{\left | \left | v_1 + \sum_{i>1}\frac{C_i}{C_1} \left [ \frac{1 - \lambda_i \Delta t}{1 - \lambda_1 \Delta t}\right ]^n e^{(\lambda_1 - \lambda_i)t_0}v_i \right | \right |_1 }\\
\end{eqnarray*}

This last expression is valid provided that $C_1 \neq 0$ and $0< 1 - \lambda_1 \Delta t.$  The first condition will be satisfied if we take $f(t_0)$ to be proportional to the vector of all $1$'s.  By the Perron-Frobenius theorem, we may take $v_1$ to contain only non-negative entries.  Then
$$\langle f(t_0),\; v_1 \rangle = C_1 e^{-\lambda_1 t_0}$$
implies $C_1 \neq 0,$ since otherwise we would have $\langle f(t_0),\; v_1 \rangle = 0,$ which is impossible.
Now, if $\Delta t < \frac{1}{\lambda_{max}},$ then $0<1 - \lambda_i \Delta t$ for all $i$.  If in addition $\lambda_2 - \lambda_1 >0$,
$$0<\frac{1 - \lambda_i \Delta t}{1 - \lambda_1 \Delta t} <1\; \forall i, \mbox{ and }$$
$$\frac{f(t_0 + n\Delta t)}{||f(t_0 + n\Delta t)||_1} \rightarrow \frac{v_1}{||v_1||_1} \mbox{ as } n \rightarrow \infty.$$

%
%
Now define the walking matrix from Theorem \ref{theorem:graphmen probs} as $E=I-H_p\Delta t,$ and recall the definition $S=||EP(t-1)_{r \not \in M}||_1$. In this fashion, $H_p$ determines the walkers' transition probability matrix for moving about the graph. Then, by definition, at any time $t>0$,
\begin{align*}
P(t)_{r \not\in m}&=\left(1+\frac{W}{S}\right)EP(t-1)_{r \not \in m}\\
&=\left(1+\frac{W}{S}\right)(I-H_p\Delta t)P(t-1)_{r \not \in m}\\
&=\frac{(I-H_p\Delta t)P(t-1)_{r \not \in m}}{S},
\end{align*}
which is exactly the normalization of $(I-H_p\Delta t)P(t-1)_{r \not \in m}$. Therefore, if we take the limit as $t \to \infty$, the RHS will converge to the normalized eigenvector of the minimal eigenvalue.	
\end{proof}

\section{Quantum Approach}


Adiabatic quantum computing (AQC) is a method by which a computer uses properties of quantum dynamics to solve problems.  The basic idea behind AQC is to encode the solution to a problem at hand into the eigenvector corresponding to the lowest eigenvalue of an hermitian matrix, $H_p.$  The matrix $H_p$  is typically referred to as the problem hamiltonian while this lowest eigenvector is referred to as ground state.  Since $H_p$ is typically difficult to work with, one constructs a simpler hermitian matrix, $H_I,$ and then evolves $H_I$ to $H_p.$  This evolution is governed by Schrodinger's equation, which takes the form:

$$\hbar \frac{\partial \psi}{\partial t} = -i\left [ (1-s(t))H_I + s(t)H_p \right ] \psi (t),$$
$$ t = [0,T],\; s(t) \in [0,1] \mbox{ and }$$
$$ s(0)=0,\; s(T)=1,\; \psi(0)=\psi_0 .$$
The speed of this evolution is controlled by the schedule function, $s(t).$   A typical choice is the constant schedule:  $s(t) = \frac{t}{T}$ for $t \in [0,T],$ where $T$ is the total integration time.  If the evolution is performed sufficiently slowly \cite{born1928beweis}, then the probability of measuring the ground state of $H_p$ is very high.  For a more in depth overview, see \cite{farhi2000quantum}.

Our approach for transforming the random walk problem into a quantum framework is based on the ideas presented in \cite{mendelson2019}.
In that work, the authors used the Laplacian for the complete graph as $H_I,$ and the so-called grounded Laplacian \cite{pirani} as $H_P.$  While the authors of that work focus on the effect of the scaling parameter $\Delta t$, we will concentrate on how the time parameter $T$ relates to convergence.  In addition, our analysis will only use a constant schedule.  However, the results described in \cite{mendelson2019} suggest that Grover's schedule is likely more optimal in terms of run time.  For a detailed comparison of the impacts from using different schedules on Grover's algorithm, see \cite{Roland_2002}.



\begin{definition}\label{def:Hamiltonian System}
Let $G_{N,\alpha}$ be the graph defined in Definition \ref{def:graph} for some distinct primes $p$ and $q$ and some totative of $\alpha$ of $N$. Let $m$ be a collection of $1 \leq k \leq N-1$ vertices which we will mark on $G_{N,\alpha}$. Let $T_F$ be the total amount of time we want for our system to evolve.\\ We define the problem Hamiltonian $H_p$ to be $L=D-A$ with the rows and columns corresponding to the vertices in $m$ deleted.\\ We also define the initial Hamiltonian $H_I$ to be $\mathbb{I}-\frac{1}{N-1-k}\mathbf{1}$, where $\mathbb{I}$ is the square Identity matrix of dimension $N-1-k$ and $\mathbf{1}$ is the square matrix of all 1's of dimension $N-1-k$.\\ We define the Hamiltonian System Evolution to be:
\begin{equation*}
H(s)=(1-s)H_I+sH_p,
\end{equation*}
where $s$ is a function of time $t$, and $s(0)=0$ and $s(T_F)=1$.
\end{definition}

Combining Definition \ref{def:Hamiltonian System} with Theorem \ref{theorem:graphmen and hamiltonian} , our problem is transformed into solving the following:
$$\hbar \frac{\partial \psi}{\partial s} = -iT_FH(s)\psi (s), \quad s \in [0,1],\; \psi(0)=\psi_0, \mbox{ where } s = \frac{t}{T_F}.$$

To mimic execution on a quantum device, we take $\psi_0$ to be the eigenvector of $H_I$ corresponding to the smallest eigenvalue, with $||\psi_0||_2 = 1.$  The authors in \cite{Jarret_2016} showed that the normalized eigenvector of the minimal eigenvalue of the problem Hamiltonian is the vector of probability amplitudes corresponding to the Hamiltonian's ground state. As a result, the final probability distribution of the walkers in the quantum setting is given by the squared amplitudes of $\psi(1)$ (in a quantum setting, the actual amplitudes are unobservable).  
Thus, the probability distribution of the walkers in the classical setting correlates to, but is not the same as the probability distribution resulting from the quantum evolution.

\section{Experiments}


%
%

\begin{figure}
\begin{center}
	\begin{subfigure}[b]{0.8\textwidth}
	\includegraphics[width=\textwidth]{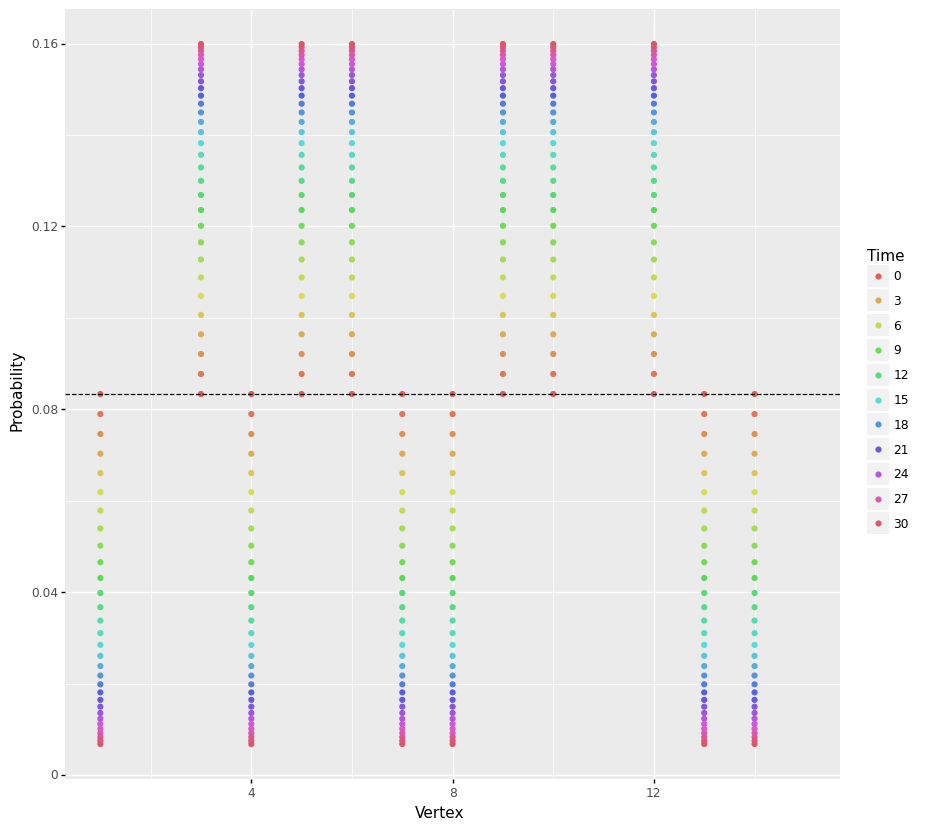}
	\caption{Classical solution for $G_{15,2} \cup G_{15,7}$}
	\label{fig:Classic 1}
	\end{subfigure}
	\begin{subfigure}[b]{0.8\textwidth}
	\includegraphics[width=\textwidth]{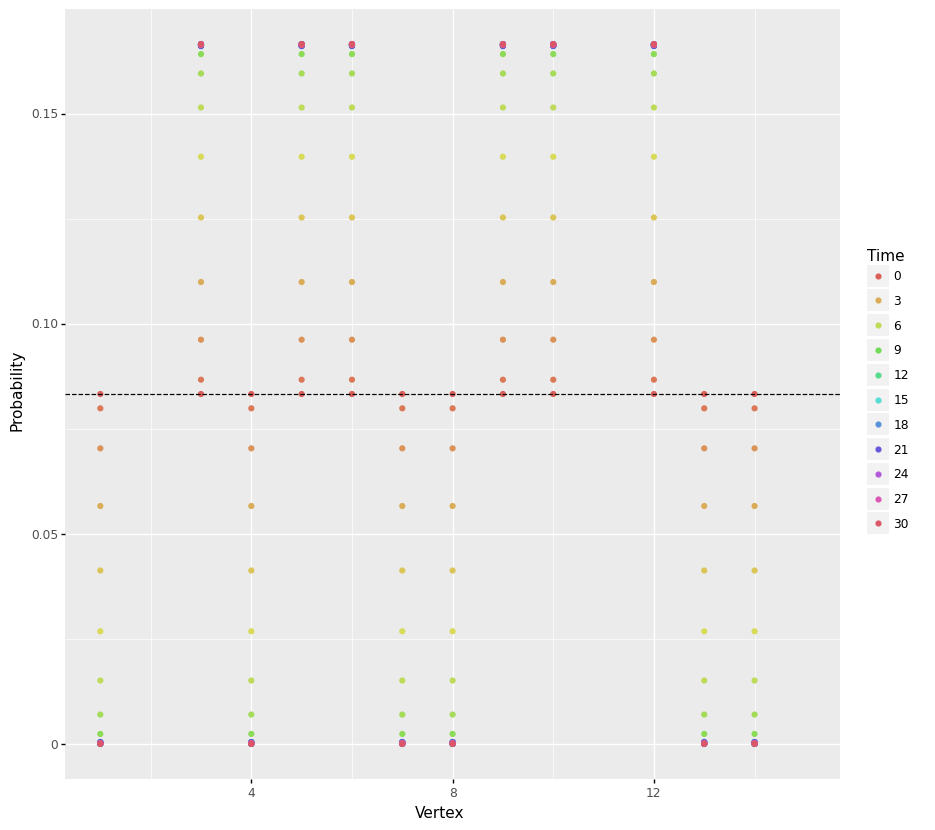}
	\caption{Quantum solution for $G_{15,2} \cup G_{15,7}$}
	\label{fig:Quantum 1}
	\end{subfigure}	
\end{center}
\caption{Classical and quantum solution for $G_{15,2} \cup G_{15,7}$}
\label{fig:compare 1}
\end{figure}

\begin{figure}
\begin{center}
	\begin{subfigure}[b]{0.8\textwidth}
	\includegraphics[width=\textwidth]{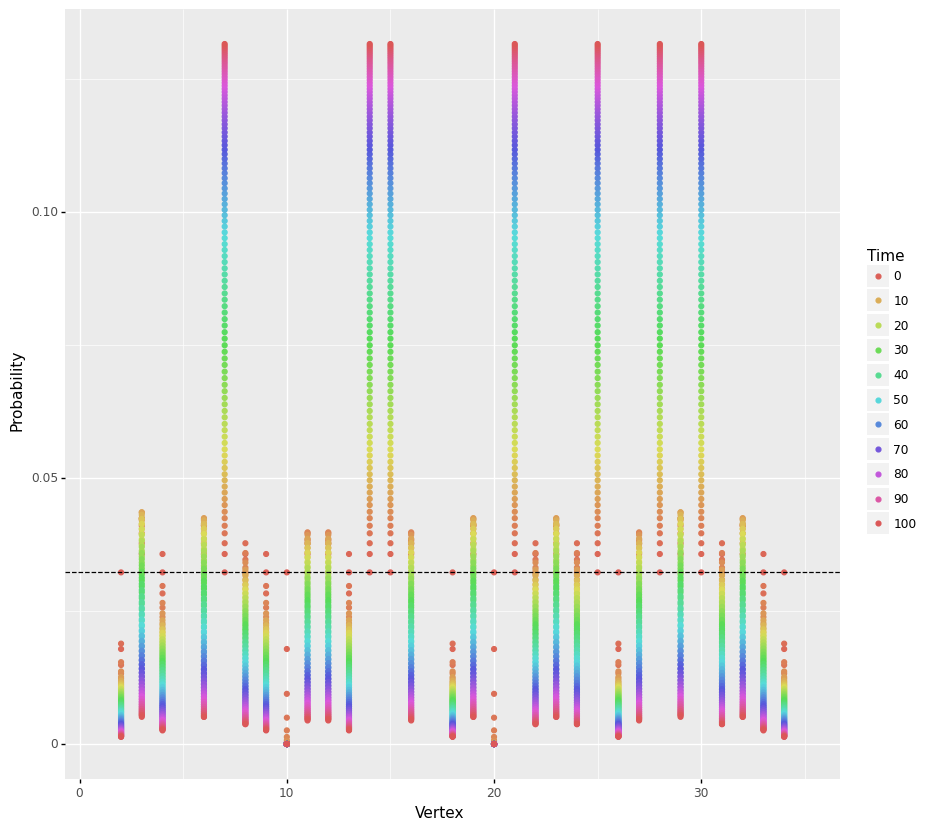}
	\caption{Classical solution for $G_{35,2}$}
	\label{fig:Classic 2}
	\end{subfigure}
	\begin{subfigure}[b]{0.8\textwidth}
	\includegraphics[width=\textwidth]{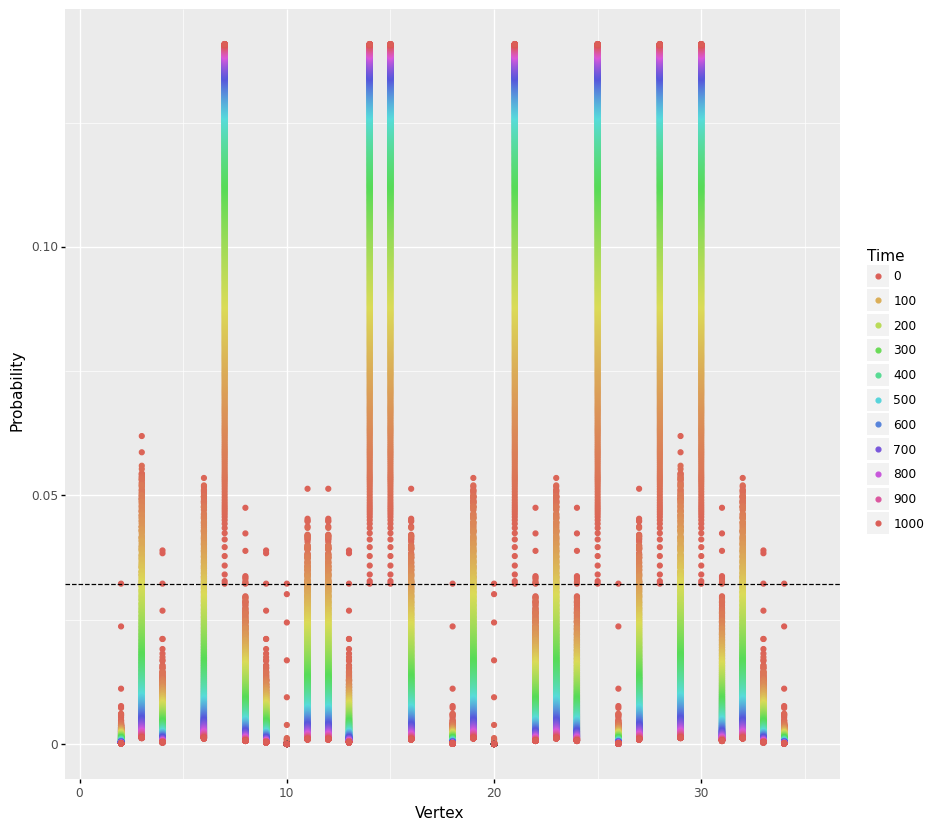}
	\caption{Quantum solution for $G_{35,2}$}
	\label{fig:Quantum 2}
	\end{subfigure}	
\end{center}
\caption{Classical and quantum solution for $G_{35,2}$.}
\label{fig:compare 3}
\end{figure}

We began our experiments by looking at the graph $G_{N,\alpha_1} \cup G_{N,\alpha_2},$ where $N = 15$ and $(\alpha_1,\alpha_2) = (2,7).$  In Figure \ref{fig:compare 1}, the dashed black line shows the initial distribution of the walkers on the graph, which is the uniform distribution about the unmarked vertices. In this case there are 12 vertices initially after selecting two vertices randomly.  We chose two vertices, 2 and 11,  based on Figure \ref{fig:2.7}, since this should maximize the probabilily of succes. The initial distribution of the walkers at each remaing vertex is $.08\overline{3}$. We see that both distributions converge to the same result, namely $.1\overline{6}$ at vertices 3,5,6,9,10, and 12, and 0 elsewhere. However, it is clear that the distribution resulting from Algorithm 2 converges more quickly. We also note that the convergence of both distributions is monotonic, meaning for each vertex, the distribution converges to the limiting distribution monotonically.  We will see later (See Figure \ref{fig:compare 3}) that this is not always the case.

%
%

\begin{figure}
\begin{center}
	\begin{subfigure}[b]{0.45\textwidth}
	\includegraphics[width=\textwidth]{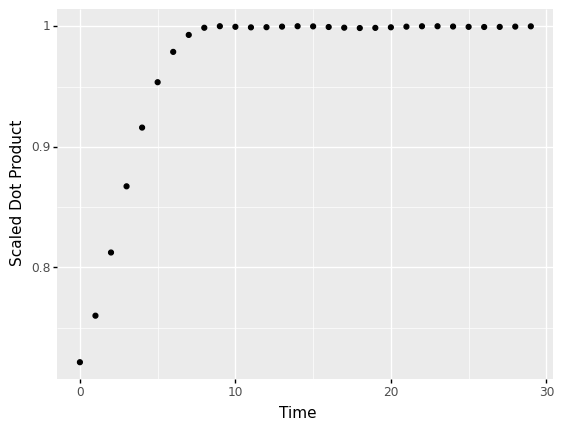}
	\caption{Scaled inner products over time for $G_{15,2} \cup G_{15,7}$ with $M = \{2,11\}$}
	\label{fig:innerproduct 1}
	\end{subfigure}
	\begin{subfigure}[b]{0.45\textwidth}
	\includegraphics[width=\textwidth]{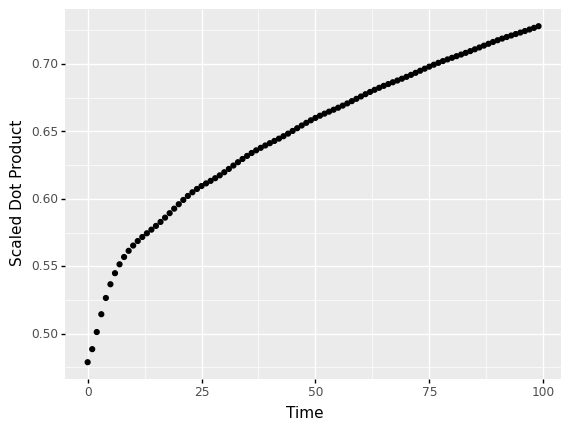}
	\caption{Scaled inner products over time for $G_{35,2}$ with $M = \{1,5,17\}$}
	\label{fig:innerproduct 2}
	\end{subfigure}	
\end{center}
\caption{Scaled inner products for different graphs using Algorithms 1 and 2}
\label{fig:compare 2}
\end{figure}

Next, we compare the distributions resulting from Algorithms 1 and 2 using the same graph and marks after sufficiently many iterations for convergence. We normalize both probability distributions, then find the cosine of the angle between them.
From Figure \ref{fig:compare 2}, we see that the distribution from Algorithm 2 converges very quickly to the limiting distribution, and for $T \geq 10$, the two distributions are practically the same.

%
%

The next graph we examine is $G_{35,2}$.  The results comparing Algorithms 1 and 2 are shown in Figure \ref{fig:compare 3}.

The dashed black line shows the initial uniform distribution of the walkers on the graph.  In this case there are 31 vertices initially, because of the 34 possible vertices, vertices 1,5, and 17 are marked. Vertices 1 and 17 were chosen randomly while vertex 5 was not.  The goal was to observe the performance in the event one of the marked vertices was a factor of $N=35.$  The initial distribution of the walkers at each remaing vertex is approximately $.0323$. We see that both distributions converge to the same result:  $.\overline{142857}$ at vertices 7,14,15,21,25,28 and 30, and 0 elsewhere. Note, however that the converence of both distributions is not monotonic.  In particular, Figure \ref{fig:compare 3} shows that the probability for several vertices (e.g., 3,4,6, and 8) increases at first and then decreases. We note that there are no vertices whose distribution first decreases then increases. We also note that the maximum value of the vertices whose distributions first increase then decrese, is greater in the Quantum Algorithm 2, than in the Classical Algorithm 1.

%
%

Figure \ref{fig:compare 3} compares the distributions from Algorithms 1 and 2 for the graph $G_{35,2}$ after sufficiently many iterations to converge.  In this case, we see that the distribution converges very slowly, especially compared to the first example which took until $T\approx 10$.  This example has not converged for $T<100$, and we can see from the Quantum Algorithm 2 Distribution in the previous image, the distribution does not appear to converge until around $T=1000$.  The reason for the slow convergence in both cases is that the spectral gap is very small ($\sim \mathcal{O}(10^{-15})$).  From the proof of Theorem \ref{theorem:graphmen and hamiltonian}, we know this will have an adverse effect on the speed of convergence.  Since the evolution time, $T,$ scales proportionally to the inverse of the spectral gap, the quantum algorithm also takes longer to converge.

\begin{figure}[h!]
\begin{center}
\includegraphics[width=14 cm, height=10 cm]{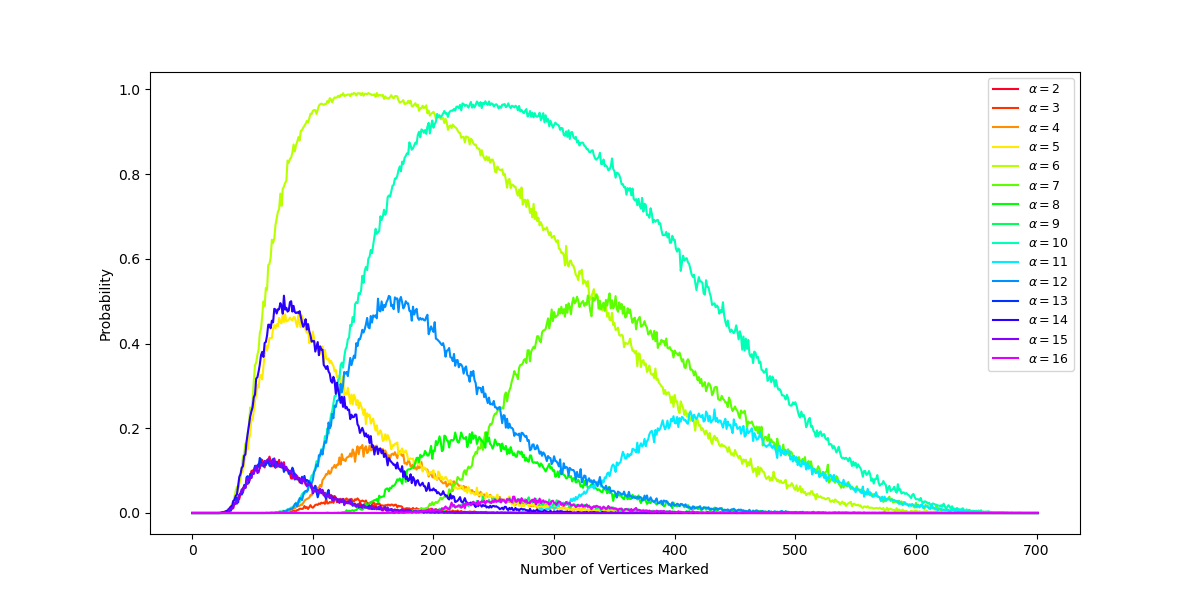}
\end{center}
\caption{Probability of satisfying conditions in Table \ref{table:conditions} for $N = 703$ and various $\alpha$'s.}
\label{fig:703}
\end{figure}

\begin{figure}
\begin{center}
	\begin{subfigure}[b]{0.9\textwidth}
	\includegraphics[width=\textwidth]{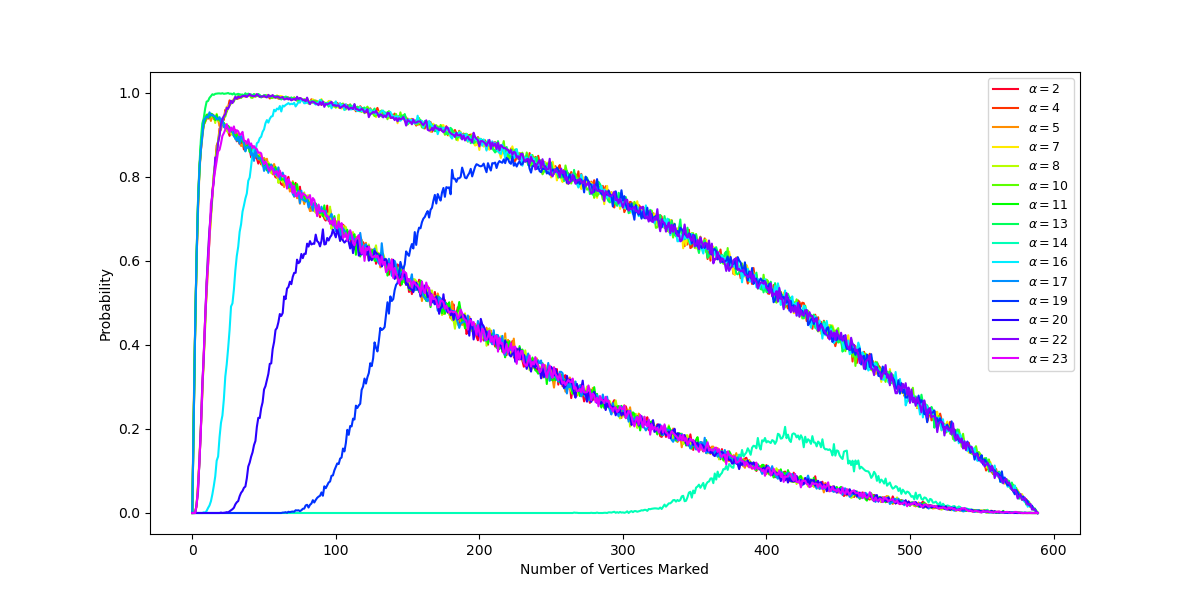}
	\caption{Probabilities for $N=591.$}
	\label{fig:591}
	\end{subfigure}
	\begin{subfigure}[b]{0.9\textwidth}
	\includegraphics[width=\textwidth]{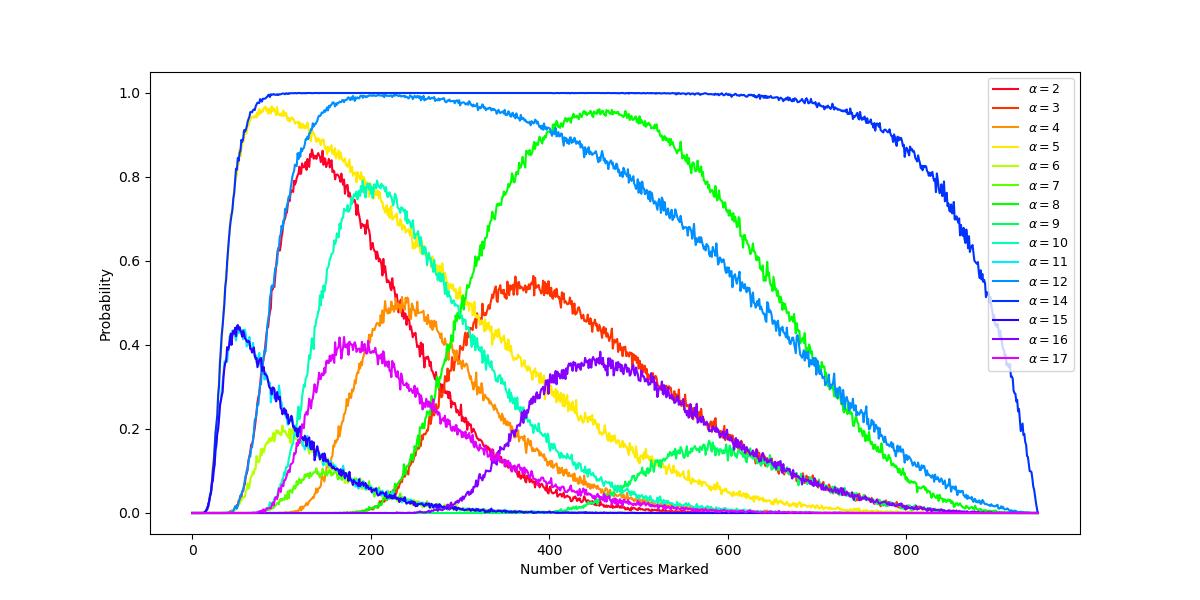}
	\caption{Probabilities for $N = 949.$}
	\label{fig:949}
	\end{subfigure}	
\end{center}
\caption{Probability of satisfying conditions in Table \ref{table:conditions} for $N = 591,949$ and various $\alpha$'s.}
\label{fig:alphas}
\end{figure}

For our next set of experiments, we select a random semiprime between 500 and 1000, and estimate the probability of satisfying the conditions in Table \ref{table:conditions} for various $\alpha$ values.   Figure \ref{fig:703} shows the results with $703 = 19 \cdot 37$. The first 13 totatives of 703 were used as $\alpha$ individually; they are listed on the right, as are their associated color. 
For each $\alpha$, 1000 random samples of $k$ marks, for $1 \leq k \leq 703$, were used. If a set of marks was successful according to the conditions in Table \ref{table:conditions}, the set was given a value of 1, and 0 otherwise. The probabilities plotted are the average success rate for each number of marks.\\
We see that there is extreme variability in not only the maximum probability of success, but also the range of highly successful marks, and distribution of the ranges of highly succesful marks for each $\alpha$. Unsurprisingly, there is significant variability between the plots for different values of $N$ (see Figure \ref{fig:alphas}).  We note that, in most cases, there appear to be several different "levels" of $\alpha$'s governed by their maximal success rate $R$, those $\alpha$'s with $0 \leq R \leq .1$, those with $.1 \leq R \leq .25$, those with $.375 \leq R \leq .625$, and those with $.875 \leq R \leq 1$. We intend to look into whether these levels persist for other semiprimes, and the other properties of the distributions of the success rate listed above, to determine whether an ``optimal" $\alpha$ can be determined for a given semiprime.

%

\begin{algorithm}
\caption{Walkers Walking Simulation}
\label{GWS}
\algsetup{indent=2em}
\begin{algorithmic}[1]
\STATE \textbf{GWS}$(N$,$a$,$m$,$\Delta t$,its$)$\\
\textbf{Input}: $N$ - semiprime, $a$ - list of $\alpha$'s, $m$ - list of vertices to mark, $\Delta t$ - scaling factor, its - number of iterations to run through.
\begin{ALC@g}
\STATE $G=$ Empty Graph
\FOR {$\alpha$ in $a$}
\STATE $G=G \cup G_{N,\alpha}$
\ENDFOR
\STATE $L=$ Laplacian$(G)$
\FOR {$v$ in $m$}
\STATE remove row $v$ and column $v$ from $L$
\ENDFOR
\STATE $E=I-L\Delta t$
\STATE $P=\frac{1}{N-1-|m|}\mathbf{1}$
\FOR {$t=1$ \TO its}
\STATE $S=||EP||_1$
\STATE $W=1-S$
\STATE $P=(1+\frac{W}{S})EP$
\ENDFOR
\end{ALC@g}
\RETURN $P$
\end{algorithmic}	
\end{algorithm}

\newpage
%

\begin{algorithm}
\caption{Quantum Amplitudes}
\label{Quantum Amplitudes}
\algsetup{indent=2em}
\begin{algorithmic}[1]
\STATE \textbf{QA}$(N$,$a$,$m$,$T)$\\
\textbf{Input}: $N$ - semiprime, $a$ - list of $\alpha$'s, $m$ - list of vertices to mark, $T$ - amount of time to run through.
\begin{ALC@g}
\STATE $G=$ Empty Graph
\FOR {$\alpha$ in $a$}
\STATE $G=G \cup G_{N,\alpha}$
\ENDFOR
\STATE $L=$ Laplacian$(G)$
\STATE $P=\mathbf{1}$
\FOR {$v$ in $m$}
\STATE replace row $v$ and column $v$ with 0's and $L[v,v]=1$
\STATE $P[v]=0$
\ENDFOR
\STATE $P=P/||P||_2$
\STATE $Probs=$ Simulate-AQC($H_F=L$, $Time=T$, Initial Condition $=P$)
\STATE $Amps=\sqrt{Probs}/||\sqrt{Probs}||_1$
\end{ALC@g}
\RETURN $Amps$
\end{algorithmic}	
\end{algorithm}

\clearpage

\section{Conclusion}
In this work we present a graph-theoretic approach to factoring semiprimes.  Beginning with a semiprime $N$, we build a graph with vertices representing the positive integers $1, \ldots, N-1.$  Edges are placed between vertices/integers that satisfy a chosen equivalence relation.  From the algebraic properties of this equivalence relation and the resulting graph, we show that one can probabilistically determine the prime factors of $N.$  To make the solution to this problem tractable, we convert the problem so that it can be solved via random walk and an equivalent adiabatic quantum approach.  The main restriction with both approaches is that the combinatorial Laplacian of the graph must possess a nonzero spectral gap.  We verify our algorithms with a series of experiments designed to not only test performance, but to compare the classical algorithm with its quantum counterpart.  We also highlight how the algorithms behave as the number of marked components increases.  It is our hope that the methods we developed, and the process itself, will be beneficial for the development of new quantum algorithms.


\bibliographystyle{plain}
\bibliography{crawbib}


\end{document}